\documentclass[12pt,oneside]{amsart}
\usepackage{amsmath,amsthm,amssymb,amscd}
\usepackage[letterpaper,left=2.5cm,right=2.5cm,top=3cm]{geometry}
\renewcommand{\epsilon}{\varepsilon}
\renewcommand{\phi}{\varphi}
\newcommand{\bZ}{\mathbb{Z}}
\newcommand{\bR}{\mathbb{R}} \newcommand{\bN}{\mathbb{N}}

\newcommand{\cT}{\mathcal{T}}
\newcommand{\cA}{\mathcal{A}}
\newcommand{\bT}{\mathbb{T}}

\newcommand{\cB}{\mathcal{B}} 
\newcommand{\brac}[1]{\langle #1\rangle}
\newtheorem{theorem}{Theorem}[section]
\newtheorem{lemma}[theorem]{Lemma}
\newtheorem{proposition}[theorem]{Proposition}
\newtheorem{prop}[theorem]{Proposition}
\newtheorem{corollary}[theorem]{Corollary}
\theoremstyle{definition}
\newtheorem{definition}[theorem]{Definition}
\newtheorem{example}[theorem]{Example}
\theoremstyle{remark}
\newtheorem{remark}[theorem]{Remark}
\numberwithin{equation}{section}

\begin{document}
\title{Trigonometric quasi-greedy bases for $L^p(\bT;w)$}
\date{\today}
\author[M.\ Nielsen]{Morten Nielsen}
\address{Department of Mathematics\\ Washington University
  \\ Campus Box 1146\\ St.\ Louis, MO~63130}
\email{mnielsen@math.wustl.edu}

\begin{abstract}
  We give a complete characterization of $2\pi$-periodic weights $w$
  for which the usual trigonometric system forms a quasi-greedy basis
  for $L^p(\bT;w)$, i.e., bases for which simple thresholding
  approximants converge in norm. The characterization implies that
  this can happen only for $p=2$ and whenever the system forms a
  quasi-greedy basis, the basis must actually be a Riesz basis.
\end{abstract}
\maketitle
\section{Introduction}
Let $\mathcal{B}=\{e_n\}_{n\in\bN}$ be a bounded Schauder basis for a
Banach space $X$, i.e., a basis for which
$0<\inf_n\|e_n\|_{X}\leq\sup_n\|e_n\|_X<\infty$. An approximation
algorithm associated with $\cB$ is a sequence $\{A_n\}_{n=1}^\infty$
of (possibly nonlinear) maps $A_n:X\rightarrow X$ such that for $x\in
X$, $A_n(x)$ is a linear combination of at most $n$ elements from
$\cB$. We say that the algorithm is convergent if
$\lim_{n\rightarrow\infty} \|x-A_n(x)\|_X=0$ for every $x\in X$. For a
Schauder basis there is a natural convergent approximation
algorithm. Suppose the dual system to $\cB$ is given by
$\{e_k^*\}_{k\in\bN}$. Then the linear approximation algorithm is
given by the partial sums $S_n(x)=\sum_{k=1}^n e_k^*(x)e_k$.

Another quite natural approximation algorithm is the greedy
approximation algorithm where the partial sums are obtained by
thresholding the expansion coefficients.  Greedy approximation
algorithms are often applied successfully in applications such as
denoising and compression using wavelets, see e.g.\
\cite{MR1162221,MR1175690}.  The algorithm is defined as follows.  For
each element $x\in X$ we define the greedy ordering of the
coefficients as the map $\rho:\bN\rightarrow \bN$ with
$\rho(\bN)\supseteq \{j: e^*_j(x)\not=0\}$ such that for $j<k$ we have
either $|\,e^*_{\rho(k)}(x)|<|\,e^*_{\rho(j)}(x)|$ or
$|\,e^*_{\rho(k)}(x)|=|\,e^*_{\rho(j)}(x)|$ and $\rho(k)>\rho(j)$.
Then the greedy $m$-term approximant to $x$ is given by
$\mathcal{G}_m(x)=\sum_{j=1}^m e^*_{\rho(j)}(x)e_{\rho(j)}$.  The
question is whether the greedy algorithm is convergent.  This is
clearly the case for an {\em unconditional basis} where the expansion
$x=\sum_{k=1}^\infty e_k^*(x)e_k$ converges regardless of the
ordering. However, Temlyakov and Konyagin \cite{MR1716087} showed that
the greedy algorithm may also converge for certain conditional
bases. This lead them to define so-called quasi-greedy bases, see
\cite{MR1716087}. The definition of a quasi-greedy basis in
\cite{MR1716087} is slightly technical, but it was shown by
Wojtaszczyk \cite{MR2001k:46017} to be equivalent to to the following
statement which we use as definition.

\begin{definition}
  A bounded Schauder basis for a Banach space $X$ is called {\em
    quasi-greedy} if\, $\lim_{m\rightarrow \infty}\|x-\mathcal{G}_m(x)\|_X=0$ for
  every element $x\in X$.
\end{definition}

In this note we consider the standard trigonometric system
$\cT:=\{(2\pi)^{-1/2}e^{ikx}\}_{k\in\bZ}$ on $\bT:=[-\pi,\pi)$. As is
very well known, $\cT$ is an unconditional (orthonormal) basis for
$L^2(\bT)$ and it is immediate that the greedy algorithm
convergences. However, we are not so fortunate when we consider $\cT$
in $L^p(\bT)$, $p\not=2$. It was proved by Temlyakov
\cite{MR1613739,MR1646563} that $\cT$ fails to be a quasi-greedy basis
for $L^p(\bT)$, $1\leq p\leq \infty$, $p\not=2$. This negative result
was also proved independently by C{\'o}rdoba and Fern{\'a}ndez for
$1\leq p<2$, see \cite{MR1613739}.  So we have to look for spaces
other than $L^p(\bT)$ if we want to extend the positive result for
$\cT$ in $L^2(\bT)$. One possible path forward is to consider the
weighted space $$L^p(\bT;w):=\big\{f:\bT\rightarrow
\mathbb{C};\|f\|_{p,w}^p=\int_{-\pi}^\pi|f(t)|^pw(t)\,dt<\infty\big\},\quad
1<p<\infty,$$ where $w$ is a non-negative $2\pi$-periodic weight.  For
a suitable choice of weight, we can make $L^p(\bT;w)$ larger or
smaller than $L^p(\bT)$.  The dual system to $\cT$ in $L^p(\bT;w)$ for
a positive weight $w$ is (at least formally)
$$\bigg\{\frac1{\sqrt{2\pi}}\frac{e^{ikx}}{w(t)}\bigg\}_{k=1}^\infty$$ 
and the expansion
relative to this system is
$$f=\frac1{{2\pi}}\sum_{k\in\bZ} \int_{-\pi}^\pi
f(t)\overline{w(t)^{-1}e^{ikt}}
 w(t)\,dt\, e^{ikt}=\frac1{{2\pi}}\sum_{k\in\bZ}\brac{f,e^{ikt}}e^{ikt},$$
where $\brac{\cdot,\cdot}$ is the standard inner product on
$L^2(\bT)$.  Thus, the greedy algorithm for $\cT$ in $L^p(\bT;w)$
coincides with the usual greedy algorithm for the trigonometric
system. Our main result in Section \ref{s3} gives a complete
characterization of the non-negative weights $w$ on $\bT:=[-\pi,\pi)$
such that $\cT$ forms a quasi-greedy basis $L^p(\bT;w)$.  The
characterizing condition is rather restrictive: we must have $p=2$,
and for $p=2$, $\cT$ forms a quasi-greedy basis $L^2(\bT;w)$ if and
only if there exists $C>0$ such that $C^{-1}\leq w(t)\leq C$.  As a
consequence, we can conclude that $\cT$ is a quasi-greedy basis
$L^2(\bT;w)$ if and only if $\cT$ is a Riesz basis for
$L^2(\bT;w)$. This is perhaps surprising since a priori, the Riesz
basis property is much more restrictive than the quasi-greedy one.  In
Section \ref{s2} we characterize the weights $w$ such that $\cT$ is a
Schauder basis for $L^2(\bT;w)$. This characterization, and our main
result in Section \ref{s3}, is given in terms of the so-called
Muckenhoupt $A_2$-condition. Finally, we consider an application to
polynomial weights in Section \ref{s4}.

\section{Trigonometric Schauder bases for $L^p(\bT;w)$}\label{s2}
In this section we give a characterization of when the trigonometric
system form a Schauder basis for $L^p(\bT;w)$. We need to 
have a Schauder basis in order for thresholding to make sense. 
The result is a direct consequence of the celebrated
result by Hunt, Muckenhoupt, and Wheeden \cite{MR0312139}.

Let us first fix the notation. Let $e_k(t):=(2\pi)^{-1/2}e^{ikt}$ and let
$\cT=\{e_{n_k}\}_{k=1}^\infty$ be the ``natural'' ordering of the
trigonometric system given by the enumeration
$\{n_k\}_{k=1}^\infty=\{0,-1,1,-2,2,\ldots\}.$
We wish consider both the symmetric partial sum operator
$$T_N(f)=\sum_{k=-N}^N \brac{f,e_k}e_k,$$
where $\brac{\cdot,\cdot}$ is the standard inner product on $L^2(\bT)$, and the
partial sum operator
$$S_N(f)=\sum_{k=1}^N \brac{f,e_{n_k}}e_{n_k}.$$
We  need the Muckenhoupt $A_p$-condition. We use the convention that $0\cdot\infty=0$.
\begin{definition}
  A nonnegative $2\pi$-periodic function $w$ is called an
  $A_p$-weight, $1<p<\infty$, if there exists a constant $K<\infty$
  such that for every interval $I\subset \bR$,
$$\bigg(\frac1{|I|}\int_I w(t)\,dt\bigg)\bigg(\frac1{|I|}
\int_I {w(t)}^{-\frac1{p-1}}\,dt\bigg)^{p-1}\leq K.$$ The family of
all $A_p$-weights is denoted $\cA_p(\bT)$.
\end{definition}
The two trivial $A_p$-weights, $w\equiv 0$ and $w\equiv \infty$, are
not interesting from our point of view since the associated
$L^p(\bT,w)$ is either trivial or far too large to be useful. We
therefore exclude the trivial weights, and notice that all the
remaning $A_p$-weights satisfy $0<w(t)<\infty$ a.e., and one easily
verifies that $w,w^{{-\frac1{p-1}}}\in L^1(\bT)$.  The following
theorem is proved in \cite{MR0312139}.
\begin{theorem}[\cite{MR0312139}]\label{HHW}
  Let $w$ be a nonnegative $2\pi$-periodic weight  and
  consider formally $T_N:L^p(\bT;w)\rightarrow L^p(\bT;w)$, $1<p<\infty$. Let
  $\|T_N\|_{p,w}$ denote the corresponding operator norm. Then
  $\sup_N\|T_N\|_{p,w}<\infty$ if and only if $w\in \cA_p(\bT)$.
\end{theorem}

We now consider the following equivalent version, which gives a nice
characterization of when $\cT$ forms a Schauder basis for
$L^p(\bT;w)$.

\begin{proposition}\label{prop}
  Let $w$ be a nonnegative $2\pi$-periodic weight on $\bT$. Then $\cT$
  is a Schauder basis for $L^p(\bT;w)$, $1<p<\infty$, if and only if
  $w\in \cA_p(\bT)$.
\end{proposition}
\begin{proof}
  First, suppose $w\in \cA_p(\bT)$.  Then $0<w(t)<\infty$ a.e.\ and
  $\cT$ span a dense subset of $L^p(\bT;w)$. The natural bi-orthogonal
  system to $\bT$ is given by $\{w(t)^{-1}e_{n_k}\}_{k=1}^\infty$
  where we notice that $w(t)^{-1}e_{n_k}\in L^{q}(\bT;w)$,
  $1/p+1/q=1$. The partial sum operator is given by
$$S_N(f)=\sum_{k=1}^N \int_{-\pi}^\pi f(t)
\overline{w(t)^{-1}e_{n_k}(t)} w(t)\,dt\, e_{n_k}=\sum_{k=1}^N\brac{f,e_{n_k}}e_{n_k},$$
so, in particular, $S_{2N+1}=T_N$ for $N\geq 1$. Also,
$$S_{2N+2}=T_N+\langle f,e_{n_{2N+2}}\rangle e_{n_{2N+2}},$$
with $$\|\langle f,e_{n_{2N+2}}\rangle e_{n_{2N+2}}\|_{p,w}\leq
C|\langle fw^{1/p},w^{-1/p}e_{n_{2N+2}}\rangle|\leq C'\|f\|_{p,w},$$ where we
used that $w,w^{-q/p}\in L^1(\bT)$. Hence, by this observation and
Theorem \ref{HHW}, we obtain $\sup_N\|S_N\|_{p,w}<\infty$ and it
follows that $\cT$ is a Schauder basis for $L^p(\bT;w)$. Next, suppose
$\cT$ is a Schauder basis for $L^p(\bT;w)$. 
Let $\{d_k\}_{k=1}^\infty\subset L^q(\bT;w)$ denote the unique dual (bi-orthogonal) system.
We claim that $d_k=w^{-1}e_{n_k}$. To verify the claim, notice that
$$c_{j,k}:=\int_{-\pi}^\pi d_k(t)\overline{e_{n_j}(t)} w(t)\,dt=\delta_{j,k},$$
where $(c_{j,k})_j$ are nothing but the Fourier coefficients of
$d_k(t)w(t)\in L^1(\bT)$. Thus, $d_k(t)w(t)=e_{n_k}(t)$ a.e.\ In
particular, since $|d_k(t)|<\infty$ a.e., $0<w(t)<\infty$ a.e., and
$d_k(t)=w(t)^{-1}e_{n_k}(t)$.  We have
$S_N(f)=\sum_{k=1}^N\brac{f,e_{n_k}}e_{n_k}.$ The fact that $\cT$ is a
Schauder basis now
gives $$\sup_N\|T_N\|_{p,w}\leq\sup_N\|S_{N}\|_{p,w} <\infty ,$$ and
we use Theorem \ref{HHW} to conclude that $w\in \cA_p(\bT)$.
\end{proof}

\begin{remark}\label{r1}
  We can move the trigonometric Schauder basis in $L^p(\bT;w)$ to
  $L^p(\bT)$ using the isometric isomorphism 
  $U:L^p(\bT;w)\rightarrow L^p(\bT)$ defined by $U(f)={w}^{1/p}f$.
  Thus,
$$\big\{{w(t)}^{1/p}e_{n_k}(t)\big\}_{k\in \bN}\quad\text{and}\quad 
\bigg\{\frac{e_{n_k}(t)}{{w(t)}^{1/p}}\bigg\}_{k\in \bN}
$$
form a biorthogonal Schauder basis system in $L^p(\bT)$ whenever $w\in\cA_p(\bT)$. 
\end{remark}

\section{Trigonometric quasi-greedy bases for $L^p(\bT;w)$}\label{s3}
Proposition \ref{prop} tells us that $\cT$ is a Schauder basis for
$L^p(\bT;w)$ if and only if $w\in \cA_p(\bT)$. In this section we
prove the main result of this note: $\cT$ can be quasi-greedy in
$L^p(\bT;w)$ only for $p=2$, and we characterize the weights $w\in
\cA_2(\bT)$ for which $\cT$ is quasi-greedy in $L^2(\bT;w)$. First, we
needs to recall some basic property of quasi-greedy bases.
 
The first result we state is due to Wojtaszczyk \cite{MR2001k:46017},
see also \cite{MR1998906}. It shows that quasi-greedy bases are
unconditional for constant coefficients.
\begin{lemma}\label{lem:woj}
Suppose $\{b_k\}_{k\in\bN}$ is a quasi-greedy basis in a Banach space $X$. Then
there exist constants $0<c_1\leq c_2<\infty$ such that for every
choice of signs $\epsilon_k=\pm 1$ and any finite subset $A\subset
\bN$ we have
\begin{equation}\label{lemma:fund}
c_1\big\|\sum_{k\in A} b_k\big\|_X \leq \big\|\sum_{k\in A}
\epsilon_k b_k\big\|_X\leq c_2 \big\|\sum_{k\in A} b_k\big\|_X.
\end{equation}
\end{lemma}
\noindent
We can use Lemma \ref{lem:woj} together with some basic facts about
the geometry of $L^p(\bT;w)$ to prove the following result.

\begin{prop}\label{prop:main}
  Suppose that the trigonometric system
  $\mathcal{T}=\{e_{n_k}\}_{k\in\bN}$ is quasi-greedy in $L^p(\bT;w)$
  for some $1<p <\infty$. Then there exist constants $0<c_1\leq
  c_2<\infty$ such that for any
  $\epsilon=\{\epsilon_k\}_{k\in\bN}\in\{-1,1\}^{\bN}$ and any finite
  subset $A\subset\bN$,
  \begin{equation}\label{eq:u}c_1|A|^{1/2}\leq \big\|\sum_{k\in A}
    \epsilon_k e_{n_k}\big\|_{L^p(\bT;w)}\leq c_2|A|^{1/2}.
\end{equation}
\end{prop}
\begin{proof}
  First we consider the case $1< p\leq 2$. Let $r_1,r_2,\ldots$ be the
  Rademacher functions on $[0,1]$ defined by
  $r_k(t)=\text{sign}(\sin(2^k\pi t))$, and take any finite subset of
  integers $A=\{k_1,k_2,\ldots, k_N\}\subset \bN$. Put
  $D_N=\sum_{l=1}^N \epsilon_{k_n} e_{n_{k_l}}$. Using
  Lemma~\ref{lem:woj}, and the fact that $L^p(\bT;w)$ has cotype $2$
  (see e.g.\ \cite[Chap.\ 3]{MR1144277}), we obtain
\[
\big\|D_N\big\|_{L^p(\bT;w)}\asymp \int_{0}^1 \big\|\sum_{n=1}^N
r_n(u)\,e_{k_n}\big\|_{L^p(\bT;w)}\,du\geq
C\big(\sum_{n=1}^N\|e_{n_{k_l}}\|_{L^p(\bT;w)}^2\big)^{1/2} \asymp
N^{1/2}.
\]
Now suppose $2\geq p<\infty$. Then $L^p(\bT;w)$ has type
$2$ (\cite[Chap.\ 3]{MR1144277}), and using Lemma~\ref{lem:woj}, we get
the estimate
$$
\big\|D_N\big\|_{L^p(\bT;w)}\asymp \int_{0}^1 \big\|\sum_{n=1}^N
r_n(u)\,e_{n_{k_l}}\big\|_{L^p(\bT;w)}\,du\leq
C\big(\sum_{n=1}^N\|e_{n_{k_l}}\|_{L^p(\bT;w)}^2\big)^{1/2}\asymp
N^{1/2}.$$
The above estimates give
$\|D_N\|_{L^2(\bT;w)}\asymp N^{1/2}$. For $1<p<2$, we notice that 
$$\|D_N\|_{L^p(\bT;w)}\leq \|D_N\|_{L^2(\bT;w)}\asymp N^{1/2},$$ and
\eqref{eq:u}
  holds in the range $1<p\leq 2$. For $2<p<\infty$, we use
$$N^{1/2}\asymp \|D_N\|_{L^2(\bT;w)}\leq \|D_N\|_{L^p(\bT;w)}$$ to reach the conclusion.
\end{proof}

A sequence $\{b_n\}_{n\in\bN}$ in a Banach space $X$ is called
democratic if there exists $D$ such that for any finite subsets
$A,B\subset\bN$ with the same cardinality $|A|=|B|$, we have
$$\big\|\sum_{k\in A} e_k\big\|_X\leq D\big\|\sum_{k\in B} e_k\big\|_X.$$
For any democratic sequence, we can define the fundamental function
\begin{equation}\label{eq:fund}\phi(n):=\sup_{A\subset \bN:|A|\leq n}\big\|\sum_{k\in A} e_k\big\|_X.
\end{equation}
Proposition \ref{prop:main} shows that whenever $\cT$ is a
quasi-greedy basis for $L^p(\bT;w)$, $\cT$ is democratic with
fundamental function $\phi(n)\asymp n^{1/2}$.  For such bases, it is
possible to prove a strong version of the Hausdorff-Young
inequality. Let us introduce some notation.

For a sequence $\{a_n\}_{n=1}^\infty$ we denote by $\{a_n^*\}$ a
non-increasing rearrangement of the sequence $\{ |a_n | \}$. Then we
define the Lorentz norms $$\|\{a_n\}\|_{2,\infty}:= \sup_n n^{1/2}
a_n^*\quad\text{ and }\quad \|\{a_n\}\|_{2,1}:= \sum_{n=1}^\infty
n^{-1/2} a_n^*.$$

The following important theorem was proved in  \cite{MR2001k:46017}.
\begin{theorem}[\cite{MR2001k:46017}]\label{lemma:EJA}
  Let $\mathcal{B}=\{b_k\}_{k\in\bN}$ be a democratic quasi-greedy
  basis for a Banach space $X$. Suppose that the fundamental function
  \eqref{eq:fund} associated with $\cB$ satisfies $\phi(n)\asymp
  n^{1/2}$.  Then there exist constants $0 < c_1 \leq c_2 < \infty$
  such that for any coefficients $\{a_k\}$
\[
c_1 \|\{a_k\}\|_{2,\infty} \leq \big\|\sum_{k\in\bN} a_k
b_k\big\|_X \leq c_2 \|\{a_k\}\|_{2,1}.
\]
\end{theorem}

\begin{remark}
  Of special interest to us is the fact that $\|\cdot\|_{2,1}$ and
  $\|\cdot\|_{2,\infty}$ assign (approximately) the same norm to flat
  sequence. More precisely, for $\mathcal{B}=\{b_k\}_{k\in\bN}$ a
  quasi-greedy basis satisfying the hypothesis of Theorem
  \ref{lemma:EJA}, there exist $c_1,c_2>0$ such that for any
  unimodular sequence $\{a_k\}_{k\in \Lambda}$, $\Lambda\subset\bN$
  (i.e., $|a_k|=1$ for $k\in\Lambda$), we have
\begin{equation}
  \label{eq:1}
c_1|\Lambda|^{1/2}\leq\big\|\sum_{k\in\Lambda} a_kb_k\big\|_X\leq 
c_2|\Lambda|^{1/2},
\end{equation}
since
$\|\{a_k\}_{k\in\Lambda}\|_{2,1}\asymp\|\{a_k\}_{k\in\Lambda}\|_{2,\infty}\asymp|\Lambda|^{1/2}$.
The estimate \eqref{eq:1} will be used below to prove our main result,
Theorem \ref{main}.
\end{remark}

\begin{theorem}\label{main}
  Let $w$ be a nonnegative $2\pi$-periodic weight. Suppose $\cT$ is a
  quasi-greedy basis for $L^p(\bT;w)$, $1<p<\infty$. Then $p=2$, $w\in A_2$, and there
  exists a positive constant $C$ such that $C^{-1}\leq w(t)\leq C$
  a.e.
\end{theorem}

\begin{proof}
  Suppose $\cT$ is a quasi-greedy basis for $L^p(\bT;w)$. Then, in
  particular, $\cT$ is a Schauder basis for $L^p(\bT;w)$ and $w\in
  A_p$ by Proposition \ref{prop}.  Now we use the Dirichlet kernel
  $D_N:=\sum_{k=1}^N e_{n_k}$ to study $w(t)$.  For each $u\in \bT$,
  we have $e_{n_k}(t-u)=e_{n_k}(t)e_{n_k}(-u)$ with $|e_{n_k}(-u)|=1$,
  and we obtain
$$D_N(t-u)=\sum_{k=1}^N e_{n_k}(-u)e_{n_k}(t).$$
Now the estimate \eqref{eq:1} gives uniformly in $u$,
\begin{equation*}
  \label{eq:est}
c_1^2N\leq   \int_{-\pi}^\pi \big|\sum_{k=1}^N e_{n_k}(t-u)\big|^2w(t)\,dt\leq
c_2^2N,
\end{equation*}
so
\begin{equation}
  \label{eq:est2}
c_1^2\leq   \int_{-\pi}^\pi \frac1N\big|\sum_{k=1}^N e_{n_k}(t-u)\big|^2w(t)\,dt\leq
c_2^2.
\end{equation}
Notice that $\frac1N\big|\sum_{k=1}^N e_{n_k}(t-u)\big|^2$ is an
approximation to the identity at the point $u$. Thus, whenever
$u\in\bT$ is a Lebesgue point of $w\in L^1(\bT)$, we obtain
\begin{equation*}
  \label{eq:est2}
c_1^2\leq   w(u)=\lim_{N\rightarrow \infty}\int_{-\pi}^\pi \frac1N\big|
\sum_{k=1}^N e_{n_k}(t-u)\big|^2w(t)\,dt\leq
c_2^2.
\end{equation*}
We conclude that $c_1^2\leq w(t)\leq c_2^2$ a.e. Now suppose
$p\not=2$.  By Proposition \ref{prop:main},
$\|D_N\|_{L^p(\bT;w)}\asymp N^{1/2}$, and it follows from H\"older's
inequality that
$$N^{1/2}\asymp \|D_N\|_p\leq \|D_N\|_1^\theta\|D_N\|_2^{1-\theta},
\qquad 1<p<2,\quad \theta=\frac2p-1,$$
or
$$N^{1/2}\asymp \|D_N\|_2\leq \|D_N\|_1^\theta\|D_N\|_p^{1-\theta},
\qquad 2<p<\infty,\quad \theta=\frac{p-2}{2p-2}.$$ In both cases we
can conclude that $\|D_N\|_{L^1(\bT;w)}\asymp N^{1/2}$ since
$\|D_N\|_{L^2(\bT;w)}\asymp N^{1/2}$. However, this is a contradiction
since we have the well-known estimate of the Lebesgue constant for
$\cT$,
$$\|D_N\|_{L^1(\bT;w)}\asymp \|D_N\|_{L^1(\bT)}\leq C\log(N),$$
where we used $c_1^2\leq w(t)\leq c_2^2$ a.e.  Thus, $\cT$
quasi-greedy implies that $p=2$, $w\in\cA_2(\bT)$ and $c_1^2\leq
w(t)\leq c_2^2$ a.e.\end{proof}

Theorem \ref{main} shows that the class of weights $w\in \cA_2(\bT)$
such that $\cT$ is a quasi-greedy basis for $L^2(\bT;w)$ is very
restrictive. In fact, the are no {\em conditional} quasi-greedy bases
for $L^2(\bT;w)$ as the following corollary shows.

\begin{corollary}
  Let $w$ be a positive $2\pi$-periodic weight for which $\cT$ is a
  quasi-greedy basis for $L^2(\bT;w)$. Then $\cT$ is a Riesz basis for
  $L^2(\bT;w)$.
  \end{corollary}
  \begin{proof}
    Suppose $\cT$ is a quasi-greedy basis for $L^2(\bT;w)$.  According
    to Theorem \ref{main}, there exists $C>0$ such that $C^{-1}\leq
    w(t)\leq C$ a.e. Hence, for any finite sequence $\{a_k\}_k$,
$$C^{-1}\int_{-\pi}^\pi \big|\sum_k a_ke_{n_k}(t)\big|^2\,dt\leq 
\int_{-\pi}^\pi \big|\sum_k a_ke_{n_k}(t)\big|^2\,w(t)\,dt\leq C
\int_{-\pi}^\pi \big|\sum_k a_ke_{n_k}(t)\big|^2\,dt$$ In particular,
$\|\sum_k a_ke_{n_k}\|_{2,w}^2\asymp \|\{a_k\}_k\|_{\ell^2}^2,$ which
shows that $\cT$ is a Riesz basis for $L^2(\bT;w)$.
\end{proof}

\section{An application}\label{s4}
Here we consider an application for general polynomial weights of the
results obtained in the previous two sections.

\begin{proposition}\label{prop2}
  Let $P$ be a polynomial of degree $n$ with $|P(-\pi)|=|P(\pi)|$. For
  $-1/n<\mu <1/n$, $\cT$ is a Schauder basis for
  $L^2(\bT;|P|^\mu)$. For such a weight $|P|^\mu$, $\cT$ is a
  quasi-greedy (and thus Riesz) basis for $L^2(\bT;|P|^\mu)$ if and
  only if $P$ has no zeros on $\bT$.
\end{proposition}
\begin{proof} 
Stein and Ricci \cite{MR890662} proved that 
for $n\in\bN$ and $0<\mu<1/n$  there exists a uniform constant $c:=c(n,\mu)$ such that
$$\int_{-1}^1|P(t)|^{-\mu}dt\leq c\bigg(\int_{-1}^1|P(t)|dt\bigg)^{-\mu},$$
where $P$ is any polynomial of degree $n$. It follows by H\"older's inequality that
$$\int_{-1}^1|P(t)|^\mu dt\leq c'\bigg(\int_{-1}^1|P(t)|dt\bigg)^{\mu}\leq
{c'}{c}\bigg(\int_{-1}^1|P(t)|^{-\mu}dt\bigg)^{-1},$$ which together
with the fact that the class of polynomials of degree $n$ is invariant
under any dilation and translation, proves that $|P|^\mu$ is in
$\cA_2(\bT)$ for $-1/n<\mu<1/n$, provided $|P(-\pi)|=|P(\pi)|$. Thus,
for $-1/n<\mu <1/n$, $\cT$ is a Schauder basis for $L^2(\bT;|P|^\mu)$.
Obviously $|P|^\mu$ is bounded on $[-\pi,\pi]$ so $\cT$ is a
quasi-greedy (and thus a Riesz) basis for $L^2(\bT;|P|^\mu)$ if and
only if $P$ has no zeros on $\bT$.
\end{proof}

\begin{example}\label{e1}
  This is the famous example by Babenko of a conditional Schauder
  basis for $L^2(\bT)$ \cite{MR0027093}.  Using Remark \ref{r1} and
  Proposition \ref{prop2}, we see that the system
  $\{|t|^{\alpha}e_{n_k}\}_{k=1}^\infty$ forms a Schauder basis for
  $L^2(\bT)$ for $0<\alpha<1/2$ since, according to Proposition
  \ref{prop2}, $|t|^{\mu}\in \cA_2$ for $-1<\mu<1$. The basis is
  conditional since $t$ has a zero on $\bT$.
\end{example}

\end{document}